\documentclass[12pt]{amsart}
\usepackage{amsmath}
\usepackage{amsthm}
\usepackage{latexsym}
\usepackage{amssymb}

\newcounter{obsctr}

\renewcommand{\theequation}{\thesection.\arabic{equation}}

\begin{document}
\baselineskip 16pt
\def\A {{\mathcal{A}}}
\def\D {{\mathcal{D}}}
\def\R {{\mathbb{R}}}
\def\N {{\mathbb{N}}}
\def\C {{\mathbb{C}}}
\def\Z {{\mathbb{Z}}}
\def\l {\ell}
\def\ml {multline}
\def\multiline {\multline}
\def\lessim {\lesssim}
%
%
%
%
%
%
\def\phi{\varphi}
\def\epsilon{\varepsilon}
\title{Analyticity and loss of
derivatives}    
 
\author{Makhlouf Derridj}
\address{5 rue de la Juviniere, 78350 Les Loges en
Josas, FRANCE}
\email{derridj@club-internet.fr}
\author{David S. Tartakoff}
\address{Department of Mathematics, University
of Illinois at Chicago, m/c 249, 851 S.
Morgan St., Chicago IL  60607, USA}
\email{dst@uic.edu}
\date{\today}
\begin{abstract} A very recent paper of Kohn studies
hypoellipticity for a sum of squares of complex vector
fields which exhibit a large loss of derivatives. We
prove analytic hypoellipticity for this operator.
\end{abstract}
\maketitle
\pagestyle{myheadings}
\markboth{M. Derridj and D.S. Tartakoff}
{Local Analyticity for complex vector fields}
\section{Introduction}
\renewcommand{\theequation}{\thesection.\arabic{equation}}
\setcounter{equation}{0}
\setcounter{theorem}{0}
\setcounter{proposition}{0}  
\setcounter{lemma}{0}
\setcounter{corollary}{0} 
\setcounter{definition}{0}
\setcounter{remark}{0}

In \cite{K2005}, J.J. Kohn proved the hypoellipticity
of the operator 
$$P = LL^*+ (\overline{z}^kL)^*(\overline{z}^kL),
\qquad L = {\partial  \over \partial z} + i\overline{z}
{\partial \over \partial t},$$
for which there is a large loss of
derivatives - indeed in the {\it a priori} estimate
one  bounds only the Sobolev norm of order $-(k-1)/2.$

We show in this note that solutions of
$Pu=f$ with $f$ real analytic are themselves real
analytic in any open set where $f$ is. 

The {\it a priori} estimate which Kohn established
for this operator and with which we will work is 
$$\|\overline{L}v\|_0^2 +
\|\overline{z}^k{L}v\|_0^2 + \|v\|^2_{-{k-1\over 2}}
\lesssim |(Pv, v)_{L^2}| + \|v\|_{-{k\over 2}}^2,
\qquad v
\in C_0^\infty.$$ 

The estimate has two interesting parts. The first two
terms on the left exhibit maximal control in
$\overline{L}$ and $\overline{z}^kL,$ but only these
complex directions. Hence in obtaining recursive
bounds for derivatives it is essential to keep one of
these vector fields available for as long as possible.
For this, we will construct a carefully balanced
localization of high powers of $T= -2i
\partial/\partial t.$  When this becomes no long
possible (even with Ehrenpreis-type cut-off functions
and the constructions of the second author in
\cite{DT1978}, \cite{DT1980}, one can only localize a
fixed, though arbitrarily high, power of $T$), one
must accept the lack of a `good' derivative
($\overline{L}$ or $\overline{z}^kL$) and use the third
term on the left of the estimate, introduce a new
localizing function, and accept the loss of the (large
but finite) number of derivatives and start the whole
process again, but with only a fraction of the
original power of $T.$

Our first observation is that we know the analyticity
of the solution for $z$ different from $0$ from the
earlier work of the second author \cite{DT1978},
\cite{DT1980} and Treves \cite{Tr1978}. Thus, modulo
brackets with localizing functions whose derivatives
are supported in the known analytic hypoelliptic
region, we take all localizing functions independent of
$z.$

Our second observation is that it suffices to bound
derivatives measured in terms of high powers of the
vector fields $L$ and $\overline{L}$ in $L^2$ norm, by
standard arguments, and indeed estimating high powers
of ${L}$ can be reduced to bounding high
powers of 
$\overline L$ and powers of $T$ of half the order, by
repeated integration by parts.  Thus our overall scheme
will be to start with high powers (order $2p$) of $L$ or
$\overline{L},$ use integration by parts and
the {\it a priori} estimate repeatedly to reduce to
treating $T^pu$ in a slightly larger set.  

And to do this, we introduce a new special
localization of
$T^p$ adapted to this problem.

The new localization of $T^p$ may be written in the
form: 

$$(T^{p_1,p_2})_\phi = \sum_{a\leq p_1 \atop b\leq
p_2}{L^a\circ 
z^a\circ 
T^{p_1-a}\circ
\phi^{(a+b)}\circ T^{p_2-b}\circ
\overline{z}^b\circ \overline{L}^b
\over a!b!},$$
Here by $\phi^{(r)}$ we mean $(-i\partial/\partial
t)^r\phi(t)$ since near $z=0$ we have seen that we may
take the localizing function independent of $z.$ Note
that the leading term (with $a+b=0)$ is merely
$T^{p_1}\phi T^{p_2}$ which equals $T^{p_1+p_2}$
on the initial open set $\Omega_0$ where $\phi \equiv
1.$

We have the commutation relations: 

$$[L, (T^{p_1,p_2})_\phi] \equiv L\circ
(T^{p_1-1,p_2})_{\phi'},$$
$$ [\overline{L},
(T^{p_1,p_2})_\phi] \equiv 
(T^{p_1,p_2-1})_{\phi'}\circ \overline{L},
$$
$$[(T^{p_1,p_2})_\phi,z] = (T^{p_1-1,p_2})_{\phi'}\circ
z,
$$
and
$$
[(T^{p_1,p_2})_\phi,\overline{z}] =
\overline{z}\circ (T^{p_1,p_2-1})_{\phi'},
$$
\vskip.1in\noindent
where the $\equiv$ denotes modulo
$C^{p_1-p_1'+p_2-p_2'}$ terms of the form
$${L^{p_1-p_1'}\circ z^{p_1-p_1'}\circ 
T^{p_1'}\circ\phi^{(p_1'+p_2'+1)}\circ T^{p_2'}
\circ\overline{z}^{p_2-p_2'}\circ
\overline{L}^{p_2-p_2'}\over (p_1-p_1')!(p_2-p_2')!}$$
with either $p_1'=0$ or
$p_2'=0,$ i.e., terms where all free $T$
derivatives have been eliminated on one side of
$\phi$ or the other. Thus if we start with
$p_1=p_2={p/ 2},$ and iteratively apply these
commutation relations, the number of $T$ derivatives not
necessarily applied to
$\phi$ is eventually at most ${p/ 2}.$ 

So we insert first $v=(T^{{p\over 2},{p\over 2}})_\phi u$ in the {\it
a priori} inequality, then bring $(T^{{p\over 2},{p\over 2}})_\phi$
to the left of $P=-L\overline{L}-
 \overline{L}z^k\overline{z}^kL$ since $Pu$ is known
and analytic. We have, omitting for now the
`subelliptic' term,
$$\|\overline{L}(T^{{p\over 2},{p\over 2}})_\phi
u\|_0^2 +
\|\overline{z}^k{L}(T^{{p\over 2},{p\over 2}})_\phi
u\|_0^2 \lesssim |(P(T^{{p\over 2},{p\over 2}})_\phi u,
(T^{{p\over 2},{p\over 2}})_\phi u)_{L^2}|$$
$$\lesssim
|((T^{{p\over 2},{p\over 2}})_\phi Pu, (T^{{p\over 2},{p\over 2}})_\phi
u)_{L^2}| + |([P,(T^{{p\over 2},{p\over 2}})_\phi] u,
(T^{{p\over 2},{p\over 2}})_\phi u)_{L^2}|$$
and, by the above bracket relations,  
$$([P,(T^{{p\over 2},{p\over 2}})_\phi] u,
(T^{{p\over 2},{p\over 2}})_\phi u)$$ 
$$= -([L\overline{L},(T^{{p\over 2},{p\over 2}})_\phi]
u, (T^{{p\over 2},{p\over 2}})_\phi u) - ([
\overline L z^k\overline{z}^k{L},(T^{{p\over 2},{p\over
2}})_\phi] u, (T^{{p\over 2},{p\over 2}})_\phi u)$$
$$\equiv -(L(T^{{{p\over 2},{p\over 2}}-1})_{\phi'}
\overline{L}u, (T^{{p\over 2},{p\over 2}})_\phi u)
-(L (T^{{p\over 2}-1,{p\over 2}})_{\phi'}
\overline{L}u, (T^{{p\over 2},{p\over 2}})_\phi u)$$ 
$$-
((T^{{p\over 2}-1,{p\over
2}})_{\phi'}\overline{L} z^k\overline{z}^k {L}u, 
(T^{{p\over 2},{p\over 2}})_\phi u)$$
$$-\sum_{k'=1}^k 
(\overline L{z}^{k'}(T^{{p\over 2},{p\over
2}-1})_{\phi'}
{z}^{k-k'}\overline z^k{L}u, 
(T^{{p\over 2},{p\over 2}})_\phi u) $$
$$-\sum_{k'=0}^{k-1}(\overline L{z}^k
\overline z^{k'} (T^{{p\over 2}-1,{p\over 2}})_{\phi'}
\overline z^{k-k'}{L}u, 
(T^{{p\over 2},{p\over 2}})_\phi u)$$ 
$$- (\overline L{z}^k
\overline z^k{L}(T^{{p\over 2},{p\over
2}-1})_{\phi'}u,  (T^{{p\over 2},{p\over
2}})_\phi u). 
$$ 
with the same meaning for $\equiv$ as above. In every
term, no powers of $z$ or $\overline{z}$ have been
lost, though some may need to be brought to the left
of the $(T^{q_1,q_2})_{\tilde{\phi}}$ 
with again no loss
of powers of $z$ or $\overline{z}$ and a further
reduction in order, every bracket reduces the order of 
the sum of the two indices $p_1$ and
$p_2$ by one (here we started with $p_1=p_2=p/2$), pick
up one derivative on $\phi,$ and leave the vector
fields over which we have maximal control in the
estimate intact and in the correct order. Thus we may
bring either
$\overline{L}z^k$ or $L$ to the right as
 $\overline{z}^kL$ or
$\overline{L},$ and use a weighted Schwarz inequality on
the result to take maximal advantage of the {\it a
priori} inequality. Iterations of all of this
continue until there remain at most $p/2$
free $T$ derivatives (i.e., the $T$ derivatives on
at least one side of $\phi$ are all `corrected' by
good vector fields) and perhaps as many as
${p/ 2} \;L$ or $\overline{L}$ derivatives, and we
may continue further until, at worst, the remaining $L$
and
$\overline{L}$ derivatives bracket two at a time to
produce more $T$'s, one at a time. After all of this,
there will be at most
$T^{3p\over 4}$ remaining. 

It is here that the final term on the left of the
{\it a priori} inequality is used, in order to bring
the localizing function out of the norm after
creating another balanced localization of $T^{3p/4}$
with a new localizing
function of Ehrenpreis type with slightly larger
support, geared to
$3p/4$ instead of to
$p.$

This means that the entire process, which reduced the
order from $p$ to $3p/4,$ or more precisely to $3p/4 +
(k-1)/2,$ is repeated, over and over, each time
essentially reducing the order by a factor of $3/4.$
After on the order of $\log_{4/3} p$ such iterations
we are reduced to a bounded number of derivatives,
and, as in \cite{DT1978} and \cite{DT1980}, all of these
nested open sets may be chosen to fit in the one open
set $\Omega_1$ where $Pu$ is known to be analytic, and
all constants chosen independent of $p$ (but
depending on $Pu$). The fact that in those
works one full iteration reduced the order by half
played no essential role - a factor of
$3/4$ would have worked just as well. 

The final estimate, as in those works, is that for all
$\alpha$ with $|\alpha|\leq p,$
$$|D^{|\alpha|}u|_{L^\infty(\Omega_0)} \leq CC^pp^p
\sim C'{C'}^pp!$$ in $\Omega_0$ with $C$ independent of
$p,$ which proves the analyticity of the solution in
$\Omega_0.$

\end{document}